\documentclass[lettersize,romanappendices,conference]{ieeeconf}

\IEEEoverridecommandlockouts                              

\usepackage{url}
\usepackage{verbatim}
\usepackage{amsmath} 
\usepackage{amsfonts,mathrsfs}
\usepackage{amssymb}
\usepackage{color}
\usepackage{graphicx}
\usepackage{epsfig}
\usepackage{enumerate}
\usepackage{algorithm,algorithmicx}
\usepackage{algpseudocode}
\usepackage[hidelinks]{hyperref}
\usepackage{empheq}
\usepackage{bm}
\usepackage{accents}
\usepackage{cite}
\usepackage{caption}
\usepackage{subcaption}
{\it}{}
{\it}{}
{\it}{}
\newtheorem{proposition}{Proposition}{\it}{}
{\it}{}
{\it}{}
\newtheorem{remark}{Remark}{\it}{}
{\it}{}

\newcommand{\mc}{\mathcal}
\newcommand{\bb}{\mathbb}
\newcommand{\R}{\bb R}

\DeclareMathAlphabet{\mathbbmsl}{U}{bbm}{m}{sl}

\newcommand{\argmin}{\operatorname{argmin}}

\newcommand{\col}{\operatorname{col}}

\newcommand{\1}{\mathbf{1}}

\newcommand{\Rmnum}[1]{\expandafter\@slowromancap\romannumeral #1@}
\newcommand{\eod}{\ensuremath{\hfill\Box}}

\begin{document}

\title{\Large \bf Approximate solutions to the optimal flow problem of multi-area integrated electrical and gas systems}

\author{Wicak Ananduta and Sergio Grammatico
	\thanks{W. Ananduta and S. Grammatico are with the Delft Center of Systems and Control (DCSC), TU Delft, the Netherlands. E-mail addresses: \texttt{\{w.ananduta, s.grammatico\}@tudelft.nl}. }
	\thanks{This work was partially supported by the ERC under research project COSMOS (802348). }
}



\maketitle
\thispagestyle{empty}
\pagestyle{empty}
\begin{abstract}
We formulate the optimal flow problem in a multi-area integrated electrical and gas system as a mixed-integer optimization problem by approximating the non-linear gas flows with piece-wise affine functions, thus resulting in a set of mixed-integer linear constraints. For its solution, we propose a novel algorithm that consists in one stage for solving a convexified problem and a second stage for recovering a mixed-integer solution. The latter exploits the gas flow model and requires solving a linear program. We provide an optimality certificate for the computed solution and show the advantages of our algorithm with respect to the state-of-the-art method via numerical simulations. 
\end{abstract}

\section{Introduction}

  Due to its high efficiency and low carbon emission, the utilization of natural gas for electricity production currently has the fastest growing rate among fossil fuels, and in fact, it now accounts for 25\% of power generation \cite{iea22}. Differently from renewable power generators that have intermittency issues, gas-fired generators are essentially dispatchable on demand. Therefore, they are used to ensure sufficient power delivery and to {complement} renewable energy sources \cite{ackermann01distributed,pepermans05distributed}. Meanwhile, natural gas is also supplied directly to households and industries, e.g. for generating thermal energy. Therefore, to secure fuel adequacy for power generation and availability for gas consumption, one {should} consider an interdependent operation of power and gas systems \cite{wen17}. 
 In this regard, an optimal gas and power flow (OGPF) problem  concerns computing economically efficient operating points of gas production units and dispatchable power generation units, including gas-fired ones that couple power and gas networks, to meet power and gas demands while satisfying operational and physical constraints \cite{wen17}.  
 
 One of the key challenges in solving OGPF problems is dealing with static nonlinear gas flow equations, which relate the gas pressures of two {connected} nodes and the gas flow between them.  While linear approximations of power flows are acceptable for {electrical} transmission networks \cite{yang19}, gas flows are typically approximated by mixed-integer (MI) linear  or second-order cone (SOC) constraints. In particular, the works in \cite{urbina07,correa14security,zhang16,wu20decentralized} use piece-wise affine (PWA) functions to approximate gas flow equations; thus, they require binary variables to indicate the active region/piece of the PWA functions and, in turn, obtain mixed integer linear constraints. 
 On the other hand, \cite{wen17,he18,qi19,he17,liu19} follow a different approach, i.e., relaxing Weymouth gas flow equation {\cite[Eq.(15)]{he18}} into MISOC inequality constraints {through} a binary variable {that indicates} the gas flow direction.  We note that when the gas flow directions are known and fixed, the MISOC model turns into a convex SOC one \cite{wang18}. Furthermore, some attempts to improve the tightness of the solution obtained via the MISOC model have {also} been proposed. In \cite{wen17,he18}, a penalty cost on the auxiliary variable that defines the gas flow inequality is introduced, {and \cite{he18} further presents} a sequential cone programming method. 
 
When an integrated electrical and gas system (IEGS) is large and consists of multiple areas, a decentralized method to solve the corresponding OGPF problem is preferred \cite{he18,qi19,wu20decentralized}. The works in \cite{he18,qi19,jia21decentralized,wu20decentralized} opt for the  alternating direction method of multipliers (ADMM)  to design a solution algorithm. Specifically, the authors of \cite{qi19} and \cite{jia21decentralized} consider linear approximations of gas flows, resulting in convex problems and allowing for a straightforward implementation of ADMM at the cost of relatively poor gas flow approximations. Meanwhile, \cite{he18} considers the MISOC gas flow model and proposes an iterative algorithm, where, at each iteration, area-based problems with MISOC constraints are solved to update the binary decisions, and then, a convex multi-area problem with SOC and coupling constraints is solved to update the continuous decision variables. On the other hand, \cite{wu20decentralized} uses the PWA model and {applies} directly ADMM to solve a mixed-integer problem distributedly, {however} without providing convergence {guarantees.}
 
In this paper, we study the OGPF problem of a multi-area IEGS, as in \cite{he18,qi19,wu20decentralized}. 
We use {a} PWA approximation of the gas flows since we can control {its estimation} accuracy, unlike the MISOC relaxation. Furthermore, we apply the mixed logical formulation of the PWA approximation based on \cite{bemporad99} to derive a set of MI linear constraints  (Sec. \ref{sec:IEGSmodel}). 
{Our} main contribution {is to} design a two-stage algorithm to compute a solution to the OGPF problem (Sec. \ref{sec:prop_approach}). In the first stage, we convexify the OGPF problem and compute a solution to this convexified problem. We use the output of the first stage to recover a mixed-integer solution by exploiting our gas flow model. Specifically, we obtain the integer part of the decision variable by using the logical constraints defining the PWA gas flow model and then we recompute the gas pressure variables by solving a linear {program} derived from the {approximated} gas flow equations. We show that the proposed algorithm can {indeed} find an exact solution. {Moreover,} when some gas flow equations are violated, we can quantify the maximum inaccuracy. Differently from existing algorithms, e.g., those in \cite{he18,wu20decentralized}, our method does not require solving a mixed-integer optimization problem. Instead, the subproblems in the two stages are convex, allowing us to apply a distributed convex optimization method. Furthermore, to justify the gas flow model choice, we compare the performance of our algorithm with methods that use the MISOC formulation (Sec. \ref{sec:sim_res}).

\paragraph*{Notation} We denote by $\R$ {($\R_{\geq 0}$)} the set of {(non-negative)} real numbers. 
The operator $\col(\cdot)$ stacks its arguments into a column vector. The operator  $\operatorname{sgn}(\cdot)$ denotes the sign function, i.e., \vspace{-5pt}
{$$ \operatorname{sgn}(a) = \begin{cases}
	1, \quad &\text{if } a > 0, \\
	0, \quad &\text{if } a = 0, \\
	-1, \quad &\text{if } a < 0.
\end{cases}$$}

\section{Multi-area optimal gas-power flow problem}
\label{sec:IEGSmodel}

We consider the OGPF problem of a multi-area IEGS, where each area is controlled independently but is coupled with the other areas through tie-lines in the electrical power network and/or tie-pipes in the gas network. 
\vspace{-6pt}
\subsection{System model}
We first provide the model of the system, in terms of cost functions and constraints.

\paragraph*{Dispatchable generators (DGs)} 
Let $\mc I^{\mathrm{dg}}$ denote the set of DGs.  The power production, denoted by $p_i^{\mathrm{dg}} \in \bb R_{\geq 0}$, is bounded by
\vspace{-3pt}
\begin{equation}
	\begin{aligned}
		\underline{p}_{i}^{\mathrm{dg}} \leq p_{i}^{\mathrm{dg}} &\leq  \overline{p}_{i}^{\mathrm{dg}}, \ \forall i \in \mc I^{\mathrm{dg}},
	\end{aligned}
	\label{eq:p_pu_bound}	
\end{equation}
where $\underline{p}_{i}^{\mathrm{dg}} < \overline{p}_{i}^{\mathrm{dg}}$ denote the minimum and maximum generator operation capacities. We classify the DGs into the subset of gas-fueled units ($\mc I^{\mathrm{gu}}$), i.e., those that use natural gas distributed through the gas network, and that of non-gas-fueled units ($\mc I^{\mathrm{ngu}}$), i.e., $\mc I^{\mathrm{dg}}:=\mc I^{\mathrm{gu}} \cup \mc I^{\mathrm{ngu}}$ and $\mc I^{\mathrm{gu}} \cap \mc I^{\mathrm{ngu}}=\varnothing$. For each gas-fueled unit, we consider a quadratic relationship between its power production and gas consumption, denoted by $d_i^{\mathrm{gu}} \in \bb R_{\geq 0}$ \cite[Eq. (27)]{he18}, yielding the following constraints:
\vspace{-3pt}
\begin{equation}
	\begin{aligned}
		d_i^{\mathrm{gu}} &\geq
		\eta_{2,i} (p_i^{\mathrm{dg}})^2 + \eta_{1,i} p_i^{\mathrm{dg}} + \eta_{0,i}, \  \text{if } i \in \mc I^{\mathrm{gu}}, \\
		d_i^{\mathrm{gu}} &=	0, \qquad \qquad \qquad \qquad \qquad \ \ \text{if }  i \in \mc I^{\mathrm{ngu}},
	\end{aligned}
	\label{eq:conv_dgu_ppu}	%
\end{equation}%
for some constant $\eta_{2,i} > 0$ and $\eta_{1,i}, \eta_{0,i} \in \bb R$. On the other hand, the power production of the non-gas-fueled units is assumed to have a quadratic economical cost \cite{wen17,he18,liu19}; thus, we have that
\vspace{-3pt}
\begin{equation}
\hspace{-6pt}	f_{i}^{\mathrm{dg}}(p_i^{\mathrm{dg}}) \hspace{-1pt}=\hspace{-1pt} 
	\begin{cases}
		c_{2,i}^{\mathrm{dg}}(p_i^{\mathrm{dg}})^2+ c_{1,i}^{\mathrm{dg}} p_i^{\mathrm{dg}} + c_{0,i}^{\mathrm{dg}}, \  \text{if } i \hspace{-2pt} \in \mc I^{\mathrm{ngu}},\\
		0, \qquad \qquad \qquad \qquad \qquad \  \hspace{2pt} \text{if } i \hspace{-2pt} \in \mc I^{\mathrm{gu}},
	\end{cases}
	\label{eq:f_ngu}
\end{equation}
for some cost parameters $c_{2,i}^{\mathrm{dg}} > 0$ and $c_{1,i}^{\mathrm{dg}}, c_{0,i}^{\mathrm{dg}} \in \bb R$.

\paragraph*{Power network} The power generated by the DGs is used to satisfy the power demands in the  electrical network, which is represented by an undirected connected graph $\mc G^{\mathrm e} := (\mc B,\mc L)$, where $\mc B:=\{b_1,b_2,\dots, b_B \}$, with $|\mc B|=B$, denotes the set of busses (nodes) and $\mc L \subseteq \mc B \times \mc B$ denotes the set of power lines (edges). We note that assuming there exists $m$ areas, the set of busses $\mc B$ is partitioned into $m$ non-overlapping subsets, i.e., $\mc B_a$, for $a=1,2,\dots,m$, each of which represents the set of busses that belong to the same area. Therefore, $\mc L$ includes the tie lines. By considering the DC power flow approximation \cite[{Eq. (1)}]{yang19}, the power balance at each bus can be written as:
\begin{equation}
	\textstyle\sum_{j \in \mc I_i^{\mathrm{dg}}}p_j^{\mathrm {dg}} - d_i^{\mathrm e} = \sum_{j \in \mc N_i^{\mathrm e}} \tfrac{1}{X_{\{i,j\}}}(\theta_i-\theta_j), \quad \forall i \in \mc B,
	\label{eq:pow_bal}
\end{equation}
where $d_i^{\mathrm e} \in \bb R_{\geq 0}$, $\theta_i \in \bb R$, {and $X_{\{i,j\}}$}  denote the electricity demand, the voltage angle of bus $i$, {and  the reactance of line $\{i,j\} \in \mc L$,} respectively, whereas $\mc I_i^{\mathrm{dg}}$ and $\mc N_i^{\mathrm e} :=\{j\mid \{i,j\} \in \mc L \}$ denote the set of DGs connected to bus $i$ and that of neighbor busses, respectively. We also bound $\theta_i$ by
\begin{align}
	\underline{\theta}_i \leq \theta_i &\leq \overline{\theta}_i, \forall i \in \mc B, \label{eq:theta_b}
\end{align}
with $\underline{\theta}_i < \overline{\theta}_i$ being the lower and upper bounds.

\paragraph*{Gas sources}  
We denote the set of gas sources (wells) by $\mc I^{\mathrm{gs}}$ and the gas production, denoted by $g_i^{\mathrm s} \in \bb R_{\geq 0}$, is limited by the production capacity, i.e.,
\begin{equation}
		\underline{g}_{i}^{\mathrm{s}} \leq g_{i}^{\mathrm{s}} \leq  \overline{g}_{i}^{\mathrm{s}}, \ \forall i \in \mc I^{\mathrm{gs}},
		\label{eq:gs_bound}	
\end{equation}
where $\underline{g}_{i}^{\mathrm{s}} < \overline{g}_{i}^{\mathrm{s}}$ denote the minimum and maximum production. Furthermore, {we consider} a linear gas production cost \cite{wen17,he18,liu19}, i.e., for some constants $ c_{1,i}^{\mathrm{gs}}, c_{0,i}^{\mathrm{gs}} \geq 0$, 
\begin{equation}
	f_i^{\mathrm{gs}} = c_{1,i}^{\mathrm{gs}} g_i^{\mathrm{s}} + c_{0,i}^{\mathrm{gs}}, \quad \forall i \in \mc I^{\mathrm{gs}}.
	\label{eq:f_gs}
\end{equation}

 \paragraph*{Gas network}
 The gas network is represented by a directed connected graph $\mc G^{\mathrm g} := (\mc N, \mc P)$, where $\mc N:=\{n_1,n_2,\dots,n_N\}$, with $|\mc N|=N$, denotes the set of gas nodes and $\mc P \subseteq \mc N \times \mc N$ denotes the set of edges, with both the edges $(i,j),(j,i) \in \mc P$ representing the pipeline that connects nodes $i$ and $j$. Similarly to the power network, $\mc G^{\mathrm g}$ is also partitioned into $m$ non-overlapping subsets, $\mc N_a$, for $a=1,2,\dots,m$, each of which represents the set of nodes in one area, and we assume that $|\mc N_a| > 1$. The gas balance at each node $i \in \mc N$ is given by
 \begin{align}
 \sum_{j \in \mc I_i^{\mathrm{gs}}}	g_j^{\mathrm{s}} - d_i^{\mathrm{g}} - \sum_{j \in \mc I_i^{\mathrm{gu}}} d_j^{\mathrm{gu}} &= \sum_{j \in \mc N_i^{\mathrm{g}}} \phi_{(i,j)}, \quad \forall i \in \mc N, \label{eq:gbal}
 \end{align}
where $\mc I_i^{\mathrm{gs}}$, $\mc I_i^{\mathrm{gu}}$, and $\mc N_i^{\mathrm{g}}:=\{j \mid (i,j) \in \mc P \}$ denote the set of gas sources located at gas node $i$, that of gas-fueled DGs connected to gas node $i$, and that of neighbors of node $i$. Moreover, $d_i^{\mathrm{g}}$ and $\phi_{(i,j)}$ denote the gas demand of node $i$ and the gas flow between nodes $i$ and $j$ observed by node $i$, {respectively.} In the literature, passive gas flows, assumed to be in the internal pipelines of each area, are typically modeled based on Weymouth equation whereas the gas flows in the tie pipelines can be actively controlled, as in \cite{he18}. Thus, we have the following gas flow constraints:
\begin{align}
		\phi_{(i,j)} &= \operatorname{sgn}(\psi_i-\psi_j) c_{(i,j)}^{\mathrm f}\sqrt{|\psi_i-\psi_j|}, \  \forall (i,j) \in \mc P^{\mathrm {nt}}, \label{eq:gf_eq} \\
		0&=\phi_{(i,j)} + \phi_{(j,i)}, \qquad \qquad \qquad \quad  \forall (i,j) \in \mc P^{\mathrm t},  \label{eq:gf_recip}
\end{align}
where $\psi_i$ denotes the gas pressure at node $i \in \mc N$ {and $c_{(i,j)}^{\mathrm f}$ denotes Weymouth constant that depends on the pipeline characteristics.} The set $\mc P^{\mathrm t} \subset \mc P$ denotes the set of tie pipelines that connect two neighboring areas while $\mc P^{\mathrm {nt}}:= \mc P\backslash \mc P^{\mathrm t}$ denotes the remaining (internal, non-tie) pipelines. In addition, we also constrain {the tie pipeline gas flows and} gas pressures as follows:
\begin{align}
	-\overline \phi_{(i,j)} \leq \phi_{(i,j)} \leq \overline \phi_{(i,j)}, \quad &\forall (i,j) \in \mc P^{\mathrm t},  \label{eq:phi_lim}\\
	\underline \psi_i  \leq \psi_i \leq \overline \psi_i,\quad  & \forall i \in \mc N, \label{eq:pres_lim} 
\end{align}
where $\overline \phi_{(i,j)}$ denotes the maximum gas flow and $\underline \psi_i < \overline \psi_i$ denote the minimum and maximum gas pressures of node $i$. 

The gas flow equations in \eqref{eq:gf_eq} are {nonlinear,} {implying non-convexity of the problem. Here, we approximate \eqref{eq:gf_eq} with $r$ pieces of affine functions, represented by a set of mixed-integer linear constraints \cite{bemporad99}, as follows:} 
\begin{align}
	h_{(i,j)}(y_{(i,j)},z_{(i,j)}  ) &= 0, \ \ \forall (i,j) \in \mc P^{\mathrm {nt}}, \label{eq:hf_cons1}\\
	g_{(i,j)}(y_{(i,j)},z_{(i,j)} ) &\leq 0, \ \ \forall (i,j) \in \mc P^{\mathrm {nt}}, \label{eq:gf_cons1}
\end{align}
where $h_{(i,j)}$ and $g_{(i,j)}$ are affine. We define  $y_{(i,j)}\hspace{-2pt}:=\hspace{-2pt} \col(\psi_i, \psi_j, \phi_{(i,j)}, y_{(i,j)}^{\psi_i}, \{y_{(i,j)}^m \}_{m=1}^r )\hspace{-2pt}  \in\hspace{-2pt} \bb R^{4+r}$, where $y_{(i,j)}^{\psi_i}$ and $ y_{(i,j)}^m$, for $m=1,\dots,r$, denote continuous extra variables, while $z_{(i,j)} := \col(\delta_{(i,j)}^{\psi_i}, \{\alpha_{(i,j)}^m, \beta_{(i,j)}^m, \delta_{(i,j)}^m\}_{m=1}^r  ) \in \{0,1\}^{1+3r}$ collects the binary decision variables. {For ease of presentation, we show} {the complete model in Appendix \ref{ap:pwa_gf}.} 

\subsection{Optimization problem formulation}
We can now state the overall optimization problem of the system. Let us first denote by $\bm u$ the collection of all decision variables, i.e., 
	$\bm u 
	:= \col(\{p_i^{\mathrm{dg}}, d_i^{\mathrm{gu}}\}_{i \in \mc I^{\mathrm{dg}}}, \{\theta_i\}_{i \in \mc B}, \{g_i^{\mathrm{s}}\}_{i \in \mc I^{\mathrm{gs}}},  
	\{\phi_{(i,j)}\}_{(i,j) \in \mc P^{\mathrm t}},$ $ \{y_{(i,j)}, z_{(i,j)}\}_{(i,j) \in \mc P^{\mathrm {nt}} } \} ).$ 
Then, we can  write a mixed-integer OGPF problem of a multi-area IEGS as follows
	\begin{subequations}
		\begin{empheq}[left={ \empheqlbrace\,}]{align}
	\underset{\bm u}{\min} \ \ & \sum_{i \in \mc I^{\mathrm{dg}}} f_i^{\mathrm{ dg}}(p_i^{\mathrm{dg}}) + \sum_{i \in \mc I^{\mathrm{gs}}} f_i^{\mathrm{ gs}}(g_i^{\mathrm{s}}) \label{eq:cost_f}\\
	\operatorname{s.t.} \quad & z_{(i,j)} \in \{0,1\}^{1+3r}, \quad \quad \forall (i,j) \in \mc P^{\mathrm{nt}},  \label{eq:bin_const}\\
	& \text{\eqref{eq:p_pu_bound}, \eqref{eq:conv_dgu_ppu}, \eqref{eq:pow_bal}--\eqref{eq:gs_bound}, \eqref{eq:gbal},  \eqref{eq:gf_recip}--\eqref{eq:gf_cons1} hold,	} \notag
\end{empheq}
\label{eq:opt_main}%
	\end{subequations}
where the cost functions $f_i^{\mathrm{ dg}}$ and $f_i^{\mathrm{ gs}}$ are defined in \eqref{eq:f_ngu} and \eqref{eq:f_gs}, respectively. {The mixed-integer problem in} \eqref{eq:opt_main} 
can be considered as an approximated problem since the Weymouth gas flow equations are substituted with a PWA model. Furthermore, 
	the cost function in \eqref{eq:cost_f} and most constraints, i.e., those in \eqref{eq:p_pu_bound}, \eqref{eq:conv_dgu_ppu}, \eqref{eq:theta_b}, \eqref{eq:gs_bound}, \eqref{eq:gbal}, and \eqref{eq:phi_lim}--\eqref{eq:gf_cons1}  can be decomposed area-wise. The constraints that couple two neighboring areas are a subset of {the} power balance constraints in  \eqref{eq:pow_bal},  for the busses connected to the tie lines, and the flow constraints in \eqref{eq:gf_recip}. 
\begin{remark}
	We can augment  Problem \eqref{eq:opt_main} temporally by considering that each decision variable is a vector {of} dimension equal to a predefined time horizon. {Our proposed approach directly applies to this augmented problem.}
	For the ease of notation and simplicity of the exposition, {here} we keep each decision variable to be scalar. \eod  
\end{remark}
\section{Proposed two-stage method}
\label{sec:prop_approach}
A common approach to design a distributed method for a multi-agent optimization problem with coupling constraints is by using the (augmented) Lagrangian method, where coupling constraints are dualized \cite[Chap. 2]{bertsekas14constrained}. Distributed algorithms derived from this approach, such as the ADMM, typically has a global convergence guarantee when the problem is convex and under rather mild conditions.
However, to our knowledge, there is no such a guarantee when the problem is mixed-integer, such as Problem \eqref{eq:opt_main}. {Therefore,} 
we propose a two-stage method that can be implemented in a distributed fashion to solve Problem \eqref{eq:opt_main}. In the first stage, we solve a convexified problem whereas, in the second stage, we recover an approximate mixed-integer solution by exploiting the PWA gas model. 

The source of non-convexity in Problem \eqref{eq:opt_main} is the binary constraints in \eqref{eq:bin_const}. Therefore, we relax these constraints by considering their convex hulls and obtain the following {convex} problem {(with equality coupling constraints):}
	\begin{subequations}
	\begin{align}
		\underset{\bm u}{\min} \ \ & \sum_{i \in \mc I^{\mathrm{dg}}} f_i^{\mathrm{ dg}}(p_i^{\mathrm{dg}}) + \sum_{i \in \mc I^{\mathrm{gs}}} f_i^{\mathrm{ gs}}(g_i^{\mathrm{s}}) \label{eq:cost_f2}\\
		\operatorname{s.t.} \quad & z_{(i,j)} \in [0,1]^{1+3r}, \quad \quad \forall (i,j) \in \mc P^{\mathrm {nt}},  \label{eq:bin_const2}\\
		& \text{\eqref{eq:p_pu_bound}, \eqref{eq:conv_dgu_ppu}, \eqref{eq:pow_bal}--\eqref{eq:gs_bound}, \eqref{eq:gbal}, and \eqref{eq:gf_recip}--\eqref{eq:gf_cons1} hold.} \notag
	\end{align}
	\label{eq:opt_main_cvx}%
\end{subequations} 
{Now,} we can resort to distributed augmented Lagrangian algorithms, such as those in {\cite[Sect. 3]{chatzipanagiotis15},\cite[Alg. 1]{chang16a},\cite[Algs. 1 \& 2]{ananduta21},}  to solve Problem \eqref{eq:opt_main_cvx}.

Let us now suppose that an optimal solution to the convexified problem in \eqref{eq:opt_main_cvx} is obtained and denoted by ${\bm u}^\circ$. In general we cannot guarantee that the computed solution {satisfies} the binary constraints \eqref{eq:bin_const}. Therefore, now we explain how {to} recover an approximate mixed-integer solution with minimal violation on the gas flow equations.

From the PWA gas flow model, particularly \eqref{eq:phi_logic1}, \eqref{eq:phi_logic} and \eqref{eq:ub_cons1} in Appendix \ref{ap:pwa_gf}, for each $(i,j) \in \mc P^{\mathrm {nt}}$, the binary variable $z_{(i,j)}$ defines some logical implications of the flow variable $\phi_{(i,j)}$.  
Given the decision $ \phi_{(i,j)}^\circ$, for each $(i,j) \in \mc P^{\mathrm{nt}}$, we can then use these constraints  to obtain the binary decision $\tilde z_{(i,j)} := \col(\tilde \delta_{(i,j)}^{\psi_i}, \{\tilde \alpha_{(i,j)}^m, \tilde \beta_{(i,j)}^m, \tilde \delta_{(i,j)}^m\}_{r=1}^m  )$ as follows:
\begin{equation}
	\begin{aligned}
		\tilde \delta^\psi_{(i,j)}  &= \begin{cases}
			1, \quad \text{if } \phi_{ij}^\circ \leq 0, \\
			0, \quad \text{otherwise}, 
		\end{cases} \\
		\tilde 	\delta_{(i,j)}^m &= \begin{cases}
			1, \quad \text{if } \underline{\phi}_{(i,j)}^m \leq  \phi_{(i,j)}^\circ \leq \overline{\phi}_{(i,j)}^m, \\
			0, \quad \text{otherwise}, \quad \text{for } m=1,\dots, r,
		\end{cases} \\
		&-\tilde  \alpha_{(i,j)}^m +\tilde  \delta_{(i,j)}^m \leq 0, \  
		-\tilde \beta_{(i,j)}^m +\tilde  \delta_{(i,j)}^m\leq 0,  
		\\ 
		&\tilde \alpha_{(i,j)}^m +\tilde  \beta_{(i,j)}^m -\tilde  \delta_{(i,j)}^m\leq 1. 
	\end{aligned}
\label{eq:solve_binary}
\end{equation}
The (binary) variables  $\delta_{(i,j)}^m$, $m=1,\dots,r$, and $\delta^\psi_{(i,j)}$ appear in the PWA gas flow equation  \eqref{eq:gf_eq4}, restated as follows:
{\small \begin{multline} \hspace{-12pt}
\sum_{m = 1}^r \hspace{-2pt} \delta_{(i,j)}^m ( a_{(i,j)}^m {\phi}_{(i,j)} \hspace{-2pt}+\hspace{-2pt} b_{(i,j)}^m) \hspace{-2pt}=\hspace{-2pt} (2 \delta^\psi_{(i,j)}\hspace{-2pt}-\hspace{-2pt}1) \psi_{i} \hspace{-2pt}-\hspace{-2pt} (2\delta^\psi_{(i,j)}\hspace{-2pt}-\hspace{-2pt}1) \psi_{j}, \label{eq:gf_eq6}%
\end{multline}}%
for each  $(i,j) \in \mc P^{\mathrm{nt}}$. As parts of a solution to the first stage problem, $\bm u^\circ$, the tuple $(\phi_{(i,j)}^\circ, \psi_i^\circ, \psi_j^\circ, z_{(i,j)}^\circ)$, where $z_{(i,j)}^\circ$ is (possibly) continuous instead of binary, satisfies  \eqref{eq:gf_eq6}. However, $(\phi_{(i,j)}^\circ, \psi_i^\circ, \psi_j^\circ, \tilde z_{(i,j)})$ might not. {Thus,}   our next step is to recompute the pressure variables $\psi_i$, for all $i \in \mc N$, {while keeping the gas flow decisions as $\phi_{(i,j)}^\circ$, for all $(i,j) \in \mc P$.} 
To that end, let us  first compactly write the pressure variable $\bm \psi = \col((\psi_{i})_{i \in \mc N})$, the binary variables $\tilde{\bm{\delta}}^\psi = \col((\tilde{{\delta}}_{(i,j)}^\psi)_{(i,j)\in \mc P^{\mathrm{nt}}})$,  $\tilde{\bm{\delta}}_h^{\mathrm{pwa}} = \col((\tilde{{\delta}}_{(i,j)}^{\mathrm{pwa}})_{(i,j)\in \mc P^{\mathrm {nt}}})$, $\tilde{{\delta}}_{(i,j)}^{\mathrm{pwa}} = \col((\tilde{{\delta}}_{(i,j)}^{m}))_{m=1}^r$, and the flow variables ${\bm{\phi}}^\circ = \col(({{\phi}}_{(i,j)}^\circ)_{(i,j)\in \mc P^{\mathrm{nt}}}) $. 
We can recompute $\bm \psi$ by solving the following convex problem: 
\begin{subequations}
	\begin{empheq}[left={\tilde{\bm \psi} \hspace{-2pt}\in \hspace{-2pt}} \empheqlbrace]{align}
		\underset{\bm \psi}{\argmin}  \ &  J_\psi(\bm \psi) \hspace{-2pt}:=\hspace{-2pt} \| \bm E(\tilde{\bm{\delta}}^\psi) \bm \psi -\bm \theta({\bm{\phi}}^\circ, \tilde{\bm{\delta}}^{\mathrm{pwa}}) \|_\infty \label{eq:qp_pressure_cost}\\		
		\text{s.t.}  \  & \bm \psi \in [ \underline{\bm \psi}, \overline{\bm \psi}], \label{eq:psi_limits}
	\end{empheq}
	\label{eq:qp_pressure}%
\end{subequations} 
where the objective function $J_\psi$ is derived from the gas flow equation \eqref{eq:gf_eq6} as we aim at minimizing its error. We denote by
$\bm E(\tilde{\bm{\delta}}^\psi) :=  \col((E_{(i,j)}(\tilde \delta_{(i,j)}^\psi))_{(i,j)\in \mc P^{\mathrm{nt}}}) \in \bb R^{| \mc P^{\mathrm{nt}}| \times N}$ the transpose of the incidence matrix of  $\mc G^{\mathrm{g}}$, since for each $(i,j)$, 
\begin{equation}
	[E_{(i,j)}(\tilde \delta_{(i,j)}^\psi)]_k = 
	\begin{cases}
		(2 \tilde \delta^\psi_{(i,j)}-1), \quad &\text{if } k=i, \\
		- (2 \tilde \delta^\psi_{(i,j)}-1), \quad &\text{if } k = j,\\
		0, \quad & \text{otherwise},
	\end{cases}
	\label{eq:E_row}
\end{equation}
where $[E_{(i,j)}]_k$ denotes the $k$-th component of the (row) vector $E_{(i,j)}$, and since $\tilde \delta^\psi_{(i,j)} \in \{0,1\}$. We note that $E_{(i,j)}(\tilde \delta_{(i,j)}^\psi) \bm \psi$ equals to the concatenation of the right-hand side of the equality in \eqref{eq:gf_eq6} over {all edges} in $\mc P^{\mathrm{nt}}$.  On the other hand, $\bm \theta({\bm \phi}^\circ, \tilde{\bm \delta}^{\mathrm{pwa}}) = \col((\theta_{(i,j)}({ \phi}_{(i,j)}^\circ, \tilde{\delta}_{(i,j)}^{\mathrm{pwa}}) )_{(i,j)\in \mc P^{\mathrm{nt}}})$ with
\begin{equation}
	\theta_{(i,j)}({ \phi}_{(i,j)}^\circ, \tilde{\delta}_{(i,j)}^{\mathrm{pwa}}) = \sum_{m = 1}^r \tilde \delta_{(i,j)}^m ( a_{(i,j)}^m {\phi}_{(i,j)}^\circ + b_{(i,j)}^m),
\end{equation}
which is equal to the concatenation of the left-hand side of the equation in \eqref{eq:gf_eq6}. In addition, the constraints in \eqref{eq:psi_limits} is obtained from \eqref{eq:pres_lim}, where $\underline{\bm \psi} = \col(( \underline{\psi}_i)_{i \in \mc I})$ and $\overline{\bm \psi} = \col((\overline{\psi}_i)_{i \in \mc I})$. Problem \eqref{eq:qp_pressure} can either be solved centrally, if a coordinator exists, or distributedly by  resorting to an augmented Lagrangian method since it can be equivalently written as a linear optimization problem.
\begin{remark}
	Since each area has more than one node in the gas network and the graph $\mc G^{\mathrm g}$ is connected, every node $i$ connected to a tie pipeline $(i,j) \in \mc P^{\mathrm t}$ must have an edge with another node, say $k$, that belongs to the same area; thus, $(i,k) \in \mc P^{\mathrm{nt}}$. This fact implies that the pressure variables of all nodes in $\mc N$ are updated in the second stage. \eod 
\end{remark}

Finally, for completeness, we also update the auxiliary continuous variables introduced to define the PWA model by using their definitions {(see Appendix \ref{ap:pwa_gf}),} i.e., 
\begin{equation}
	\begin{aligned}
		\tilde{y}_{(i,j)}^{\psi_i} &= \tilde{\delta}_{(i,j)}^{\psi} \tilde{\psi}_i, \quad \quad \ \forall (i,j) \in \mc P^{\mathrm{nt}},\\
		\tilde{y}_{(i,j)}^m &= \tilde{\delta}_{(i,j)}^m \phi_{(i,j)}^\circ, \quad \hspace{2pt} \forall (i,j) \in \mc P^{\mathrm{nt}}, \ \ m=1,\dots,r.
	\end{aligned}
\label{eq:y_aux_update}
\end{equation}

	\begin{algorithm}
		\caption{Two-stage method for Problem \eqref{eq:opt_main}}
		\label{alg:prop}
		\textbf{Stage 1: Convexification}
		\begin{itemize}
			\item Compute a solution to the convexified problem in \eqref{eq:opt_main_cvx} ($\bm u^\circ$).
		\end{itemize}
		\textbf{Stage 2: Solution recovery}
		\begin{itemize}
			\item Obtain binary variable $\tilde{\bm z}$ from $\bm \phi^\circ$ via \eqref{eq:solve_binary}.
			\item Recompute pressure variable $\tilde{\bm \psi}$ by solving Problem \eqref{eq:qp_pressure}
			\item Update auxiliary variables $(\tilde{y}_{(i,j)}^{\psi_i}, \{{y}_{(i,j)}^m\}_{m=1}^r)$, for all $(i,j) \in \mc P^{\mathrm{nt}}$, according to \eqref{eq:y_aux_update}.
		\end{itemize}
		\textbf{Return:}
		\begin{equation}
			\begin{aligned}
				\bm u^\star :=& \col(\{(p_i^{\mathrm{dg}}, d_i^{\mathrm{gu}})^\circ\}_{i \in \mc I^{\mathrm{dg}}}, \{\theta_i^\circ\}_{i \in \mc B}, \{(g_i^{\mathrm{s}})^\circ\}_{i \in \mc I^{\mathrm{gs}}},   \\
				& \qquad \quad \ \{\phi_{(i,j)}^\circ\}_{(i,j) \in \mc P^{\mathrm t}}, \{y_{(i,j)}^\star, \tilde z_{(i,j)}\}_{(i,j) \in \mc P^{\mathrm {nt}} } \} ), \\
				y_{(i,j)}^\star :=& \col(\tilde \psi_i, \tilde \psi_j, \phi_{(i,j)}^\circ, \tilde y_{(i,j)}^{\psi_i}, \{\tilde y_{(i,j)}^m \}_{m=1}^r ) \  \forall (i,j) \in \mc P^{\mathrm{nt}}.
			\end{aligned}
			\label{eq:out_alg}
		\end{equation}
	\end{algorithm}

{The} proposed method {is summarized} in Algorithm \ref{alg:prop}, whose solutions 
we formally characterize next.

\begin{proposition}
	\label{prop:0opt_is_MIGNE}
		Let $\bm u^\star$ be the outcome of Algorithm \ref{alg:prop} as defined in \eqref{eq:out_alg}.  If the optimal value of the cost function in Problem \eqref{eq:qp_pressure} is zero, i.e., $J_\psi(\tilde{\bm \psi}) = 0$, then, $\bm u^\star$ is a solution to Problem \eqref{eq:opt_main}. \eod
\end{proposition}

When the optimal cost of the optimization problem in \eqref{eq:qp_pressure} is positive, this value determines the maximum violation of the gas flow PWA model. 
In addition, one can also evaluate the approximation quality with respect to the nonlinear Weymouth  equation by using the following gas flow deviation metric derived from \eqref{eq:gf_eq}: For each $(i,j) \in \mc P^{\mathrm{nt}}$,
{ \begin{equation}
	\Delta^\phi_{(i,j)} = \frac{\phi_{(i,j)}^\circ- \operatorname{sgn}(\tilde \psi_i-\tilde \psi_j) c_{(i,j)}^{\mathrm f}\sqrt{|\tilde \psi_i-\tilde \psi_j|}}{ \operatorname{sgn}(\tilde \psi_i-\tilde \psi_j) c_{(i,j)}^{\mathrm f}\sqrt{|\tilde \psi_i-\tilde \psi_j|}}.  \label{eq:gf_dev}
\end{equation}}
.

\section{Numerical simulations}
\label{sec:sim_res}

We show the performance of Algorithm \ref{alg:prop} via numerical simulations on two test cases adapted from  the (medium) 73-bus-30-node-3-area  and (large) 472-bus-40-node-4-area  networks \cite[Sect. IV]{he18}. We run 100 Monte Carlo simulations for each test case where the power and gas demands are {randomly} varied. 
All simulations\footnote{The data sets and codes for the simulations are available at https://github.com/ananduta/iegs} are carried out in Matlab {R2020b} with  Gurobi solver on a laptop computer with 2.3 GHz intel core i5 processor and 8 GB of memory.

\begin{figure}
	\centering
	\includegraphics[width=.95\columnwidth]{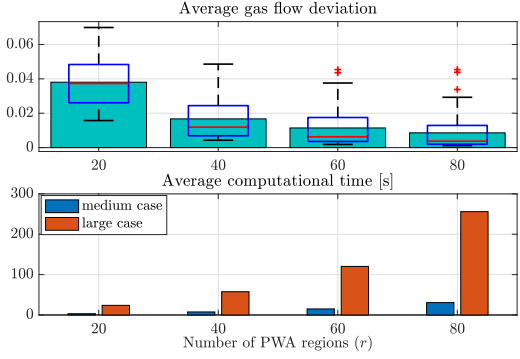}
	\caption{\small The top plot shows the average gas flow deviation, i.e.,  $\tfrac{1}{|\mc P^{\mathrm{nt}}|}\sum_{(i,j) \in \mc P^{\mathrm{nt}}}\Delta_{(i,j)}^\phi$ where $\Delta_{(i,j)}^\phi$ is defined in \eqref{eq:gf_dev}, when optimal solutions are found. 
		The bottom plot shows the average of computational time.}
	\label{fig:sim2}
\end{figure}

Figure \ref{fig:sim2} shows the performance of Algorithm \ref{alg:prop} with different numbers of PWA regions ($r$). At least in $82\%$ of the total random cases generated, the proposed algorithm with different values of $r$ finds an optimal solution to Problem \eqref{eq:opt_main}, i.e., zero cost value in the second stage ($J_\psi(\tilde{\bm \psi})=0$). In these cases, as expected, the gas flow deviations \eqref{eq:gf_dev} decrease as $r$ increases since the PWA model approximates the Weymouth equation better with larger $r$ (top plot of Figure \ref{fig:sim2}). However, when $r$ increases, consequently so does the dimension of the decision variables, the computational time also grows (bottom plot of Figure \ref{fig:sim2}).

\begin{figure}
	\hspace{-10pt}
	\begin{subfigure}{0.49\columnwidth}
		\centering
		\includegraphics[width=\textwidth]{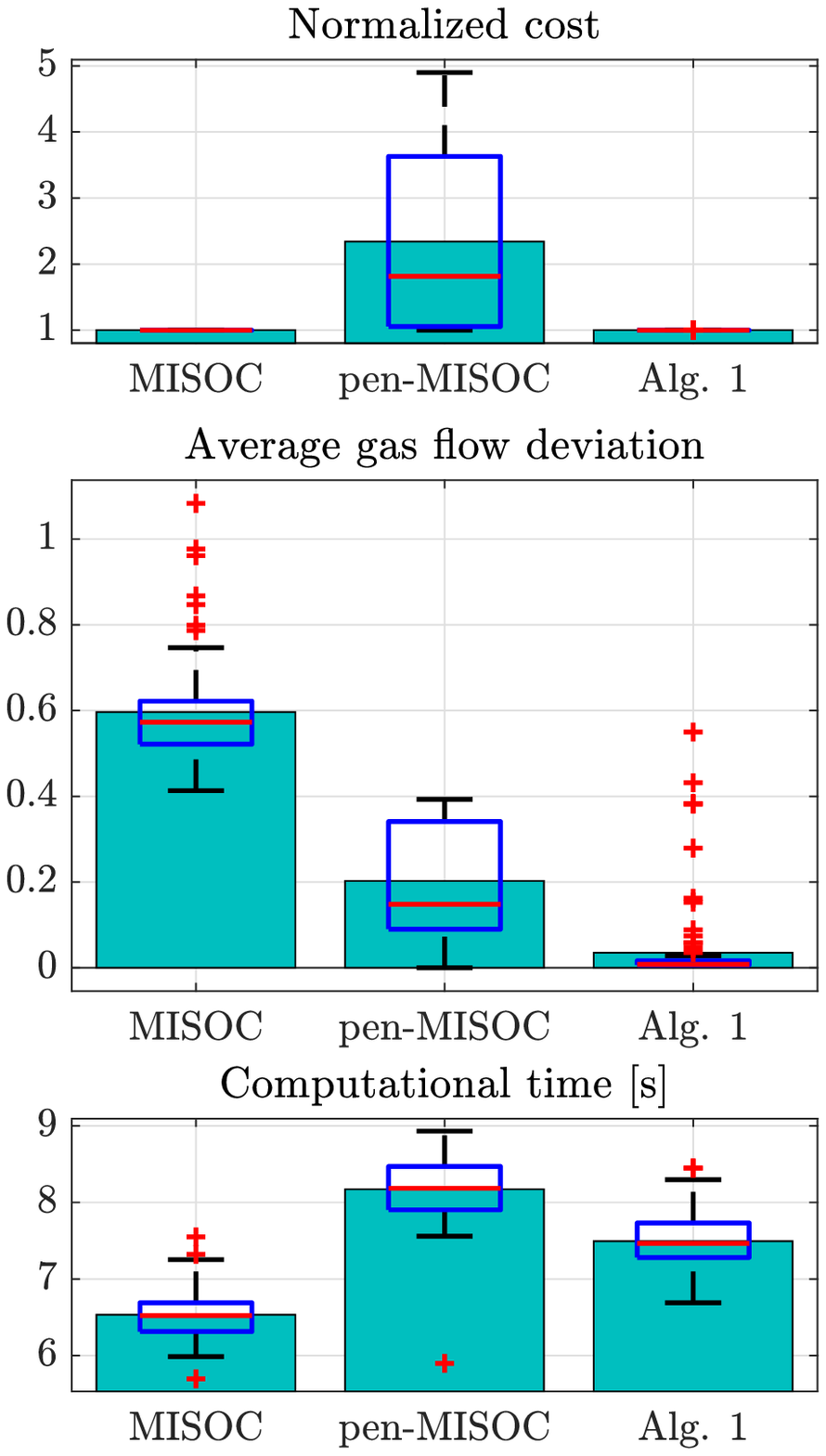}
		\caption{Medium test case}
		\label{fig:sim1_med}
	\end{subfigure}
\begin{subfigure}{0.49\columnwidth}
	\centering
	\includegraphics[width=\textwidth]{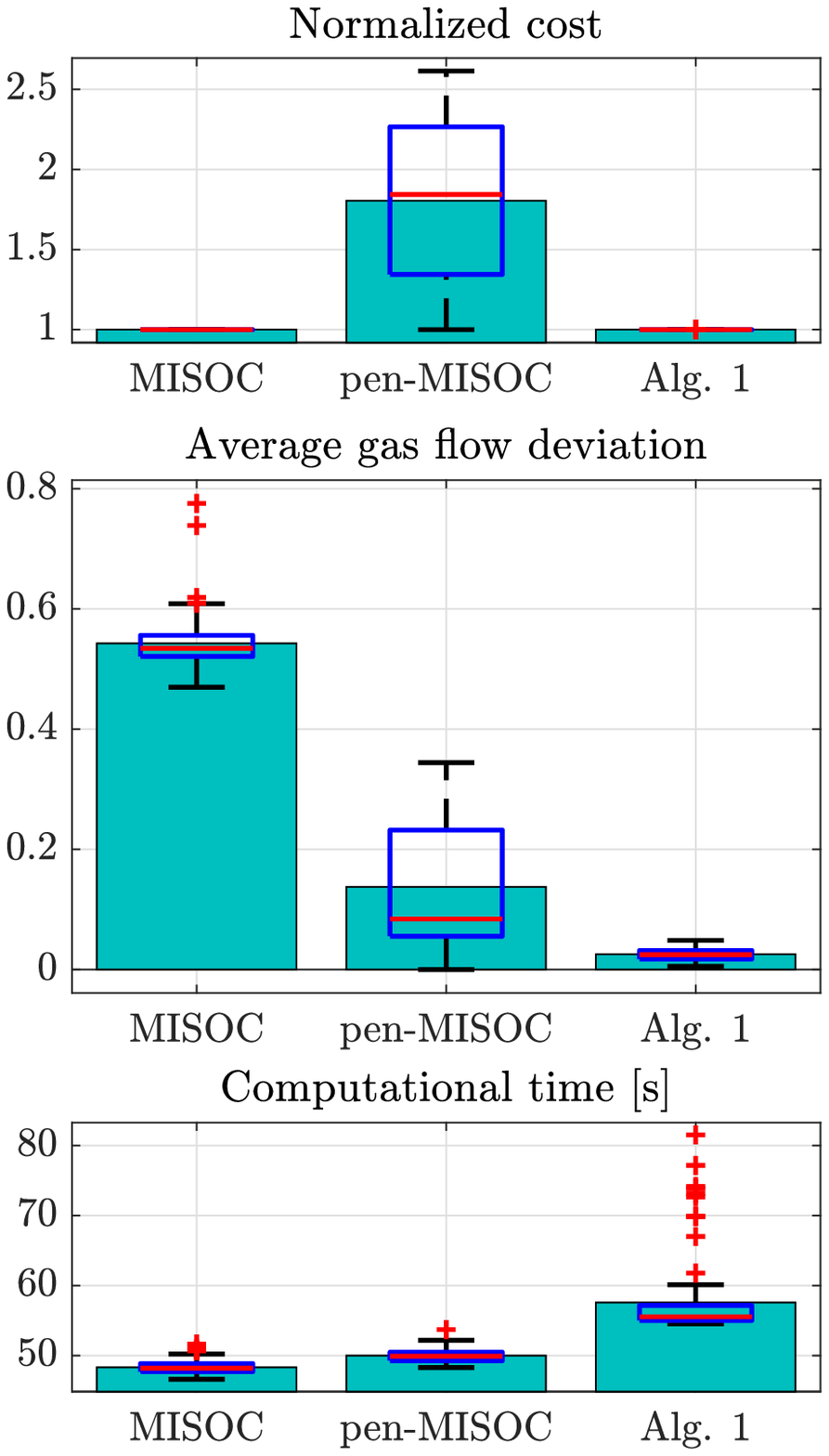}
	\caption{Large test case}
	\label{fig:sim1_large}
\end{subfigure}
\caption{\small Performance comparison between Algorithm \ref{alg:prop}, the standard MISOC, and the penalty-based MISOC (pen-MISOC) on the medium and large networks. The top plots show the (normalized) cost values, the middle plots show the average gas flow deviation, and the bottom plots show the total computational time.}
\label{fig:sim1}
\end{figure}

Then, we compare Algorithm \ref{alg:prop} with {the} state-of-the-art method based on the MISOC gas flow model \cite[Algorithm 3]{he18}, where the mixed-integer linear constraints \eqref{eq:hf_cons1}--\eqref{eq:gf_cons1} are replaced with MISOC constraints to approximate the Weymouth equations in \eqref{eq:gf_eq}.  We also consider the penalty-based variant, where an additional penalty cost function aimed at reducing gas flow deviations is introduced \cite{wen17,he18}. 
The simulation results are illustrated in Figure \ref{fig:sim1}. We can observe that Algorithm \ref{alg:prop} with $r=40$ has the smallest average gas flow deviation while achieving the same performance as the standard MISOC in terms of cost values. The penalty cost in the MISOC method can indeed reduce the average gas flow deviation although it is still larger than that of Algorithm \ref{alg:prop} while having relatively large cost values. The performance advantages of Algorithm \ref{alg:prop} comes at the cost of  total computational time, especially for the large test case.  In this regard, a trade off between the gas flow deviation and the computational time can be made by adjusting $r$, the number of regions in the PWA approximation (see Figure \ref{fig:sim2}).

\section{Conclusion}
The optimal flow problem of a multi-area integrated electrical and gas system can be formulated as a mixed-integer optimization when the nonlinear gas flow equations are approximated with piece-wise affine functions. Our  proposed algorithm can compute a solution by exploiting convexification and {the approximated gas flow model.} {Numerical} simulations show that the proposed algorithm outperforms state-of-the-art methods which {use} a mixed-integer second order cone gas flow model. Our ongoing work includes  improving the proposed algorithm, in terms of solutions and computational efficiency, and extending the problem to a generalized game setup, where selfish {yet} coupled agents exist in the network. 

\appendix

\subsection{Piece-wise affine approximation of gas flow equations}
\label{ap:pwa_gf}	 

We approximate the gas flow equality constraint at each edge $(i,j) \in \mc P^{\mathrm t}$ in  \eqref{eq:gf_eq} with a set of mixed-integer linear constraints. To that end, we introduce an auxiliary variable  $\varphi_{(i,j)} := \tfrac{ \phi_{(i,j)}^2}{	(c_{(i,j)}^{\mathrm f})^2}$ and  rewrite \eqref{eq:gf_eq} as follows:
\begin{equation}
	\varphi_{(i,j)} = \begin{cases}
		(\psi_i - \psi_j) \quad \text{if } \psi_i \geq  \psi_j, \\
	(\psi_j - \psi_i) \quad \text{otherwise.}
	\end{cases}
	\label{eq:gf_eq2}
\end{equation}
Next, we define a binary variable $\delta^\psi_{(i,j)} \in \{0,1\}$ based on the following logical constraints:
\begin{align}
	[\delta^\psi_{(i,j)}  =1] &\Leftrightarrow [ \psi_i \geq \psi_j], \label{eq:psi_logic}\\
	[\delta^\psi_{(i,j)}  =1] &\Leftrightarrow [ \phi_{ij} \geq 0].  \label{eq:phi_logic1}
\end{align}
Therefore, \eqref{eq:gf_eq2} can be rewritten as 
\begin{align}
	\varphi_{(i,j)} &= \delta^\psi_{(i,j)} (\psi_i - \psi_j) +(1-\delta^\psi_{(i,j)})(\psi_j - \psi_i )  \notag\\
	&=2 \delta^\psi_{(i,j)} \psi_i - 2\delta^\psi_{(i,j)} \psi_j + (\psi_j - \psi_i). \label{eq:gf_eq3}
\end{align} 
Let us then approximate the quadratic function $\varphi_{(i,j)} = \tfrac{ \phi_{(i,j)}^2}{	(c_{(i,j)}^{\mathrm f})^2}$ with a piece-wise affine function. Specifically, we divide the operating region of the gas flow into $r$ subregions and use a binary variable 	$\delta_{(i,j)}^m$, for each $m\in \{1,\dots,r\}$, to indicate  which subregion is active, i.e.,
\begin{align}
	[\delta_{(i,j)}^m=1] \Leftrightarrow [\underline{\phi}_{(i,j)}^m \leq \phi_{(i,j)} \leq \overline{\phi}_{(i,j)}^m], \label{eq:phi_logic}
\end{align}
with $-\overline{\phi} = \underline{\phi}_{(i,j)}^1< \overline{\phi}_{(i,j)}^1=\underline{\phi}_{(i,j)}^2< \dots <\overline{\phi}_{(i,j)}^r=\overline{\phi}$. Thus, we have the following approximation:
\begin{equation}
	\varphi_{(i,j)} \approx \sum_{m = 1}^r \delta_{(i,j)}^m ( a_{(i,j)}^m {\phi}_{(i,j)} + b_{(i,j)}^m), \label{eq:pwa_approx_varphi}
\end{equation}
for some $a_{(i,j)}^m, b_{(i,j)}^m \in \bb R$. 
Next, by \eqref{eq:gf_eq3} and \eqref{eq:pwa_approx_varphi}, we get the (approximated) gas flow equality constraint:
\begin{equation} \small 
	\sum_{m = 1}^r \delta_{(i,j)}^m ( a_{(i,j)}^m {\phi}_{(i,j)} + b_{(i,j)}^m) = 2 \delta^\psi_{(i,j)} \psi_i - 2\delta^\psi_{(i,j)} \psi_j + \psi_j - \psi_i. \label{eq:gf_eq4}
\end{equation}
Then, we introduce some auxiliary variables to substitute the products of two decision variables, i.e., $y_{(i,j)}^m:=\delta_{(i,j)}^m{\phi}_{(i,j)}$, for $m=1,\dots,r$, and  $y_{(i,j)}^{\psi_i} =\delta^\psi_{(i,j)} \psi_i$. We observe that $ \delta_{(i,j)}= 1-\delta_{(j,i)},$ and 
$\delta_{(j,i)}\psi_j = y_{(j,i)}^{\psi_j}$. 
{By combining the preceding two relationships with \eqref{eq:gf_eq4},}
 it holds that:
\begin{equation}
	\sum_{m = 1}^r  ( a_{(i,j)}^m y_{(i,j)}^m + b_{(i,j)}^m\delta_{(i,j)}^m) = 2 y_{(i,j)}^{\psi_i} + 2y_{(j,i)}^{\psi_j} - \psi_i - \psi_j, \label{eq:gf_eq5}
\end{equation}
which, together with the reciprocity constraint on $\phi_{(i,j)}$, i.e., $\phi_{(i,j)}  + \phi_{(j,i)}  = 0$, and the simplex constraint on $\delta_{(i,j)}^m$, for  all $m=1,\dots,r$, i.e., $\sum_{m=1}^r \delta_{(i,j)}^m = 1$, {is used to define} $h_{(i,j)}$ in \eqref{eq:hf_cons1}.  
Moreover, for each $(i,j) \in \mc P^{\mathrm t}$, we have:\\
1. Inequality constraints equivalent to \eqref{eq:psi_logic}:
\begin{align}
	\begin{cases}
		- \psi_i+\psi_j  \leq -(\underline \psi_i-\overline \psi_j)(1-\delta_{(i,j)}^\psi), \\
		- \psi_i+\psi_j \geq \varepsilon \1 + (-(\overline \psi_i-\underline \psi_j) - \varepsilon \1) \delta_{(i,j)}^\psi. 
	\end{cases}
\label{eq:psi1}
\end{align}
2. Inequality constraints equivalent to \eqref{eq:phi_logic1}:
\begin{align}
	\begin{cases}
			-\phi_{(i,j)} \leq \overline{\phi}_{(i,j)} (1-\delta_{(i,j)}^\psi), \\
		-\phi_{(i,j)}  \geq \varepsilon \1 + (-\overline{\phi}_{(i,j)} - \varepsilon \1) \delta_{(i,j)}^\psi. 
	\end{cases}
\label{eq:phi1}%
\end{align}
3. Inequality constraints equivalent to  \eqref{eq:phi_logic} \cite[\hspace{-1pt} (4e) \& (5a)]{bemporad99}: 
\begin{align} \small \hspace{-2pt}
	\begin{cases}
	\phi_{(i,j)} - \overline \phi_{(i,j)}^m \leq (\overline \phi_{(i,j)}-\overline \phi_{(i,j)}^m)(1-\alpha_{(i,j)}^m), \\
	\phi_{(i,j)} - \overline \phi_{(i,j)}^m \geq \varepsilon \1 + (-\overline \phi_{(i,j)} -\overline \phi_{(i,j)}^m - \varepsilon \1) \alpha_{(i,j)}^m, \\
	- \phi_{(i,j)} + \underline \phi_{(i,j)}^m \leq (\overline \phi_{(i,j)}+ \underline \phi_{(i,j)}^m)(1-\beta_{(i,j)}^m), \\
	-\phi_{(i,j)} + \underline \phi_{(i,j)}^m \geq \varepsilon \1 + (-\overline \phi_{(i,j)}+ \underline \phi_{(i,j)}^m - \varepsilon \1) \beta_{(i,j)}^m, \\
	-\alpha_{(i,j)}^m + \delta_{(i,j)}^m\leq 0, \ \ 
	-\beta_{(i,j)}^m + \delta_{(i,j)}^m\leq 0,  
	\\
	\alpha_{(i,j)}^m + \beta_{(i,j)}^m - \delta_{(i,j)}^m\leq 1, 
		\end{cases}
	\label{eq:ub_cons1}
\end{align}
for $m=1,\dots,r$, where $\alpha_{(i,j)}^m, \beta_{(i,j)}^m \in \{0,1\}$, for $m=1,\dots, r,$ are additional binary variables.\\
4. Inequality constraints equivalent to $y_{(i,j)}^m=\delta_{(i,j)}^m{\phi}_{(i,j)}$ \cite[Eq. (5b)]{bemporad99}: For all $m=1,\dots, r$,
\begin{align} \small \hspace{-8pt}
	\begin{cases}
		y_{(i,j)}^m \geq -\overline \phi_{(i,j)} \delta_{(i,j)}^m, \ \ y_{(i,j)}^m \leq \phi_{(i,j)} + \overline \phi_{(i,j)}(1-\delta_{(i,j)}^m), 
		\\
		y_{(i,j)}^m \leq \overline \phi_{(i,j)} \delta_{(i,j)}^m, \ \
		y_{(i,j)}^m \geq \phi_{(i,j)} - \overline \phi_{(i,j)}(1-\delta_{(i,j)}^m). 
	\end{cases} \hspace{-14pt}
	\label{eq:z1}
\end{align}
5. Inequality constraints equivalent to $y_{(i,j)}^{\psi_i}=\delta^\psi_{(i,j)} \psi_i$:
\begin{align}
	\begin{cases}
		y_{(i,j)}^{\psi_i} \geq \underline \psi_i \delta^\psi_{(i,j)}, \ \ y_{(i,j)}^{\psi_i} \leq \psi_i - \underline \psi_i(1-\delta^\psi_{(i,j)}), 
		\\
		y_{(i,j)}^{\psi_i} \leq \overline \psi_i \delta^\psi_{(i,j)}, \ \
		y_{(i,j)}^{\psi_i} \geq \psi_i - \overline \psi_i(1-\delta^\psi_{(i,j)}). 
	\end{cases}
	\label{eq:z3}
\end{align}
We can compactly write \eqref{eq:psi1}--\eqref{eq:z3} as $g_{(i,j)}$ in \eqref{eq:gf_cons1}.

\subsection{Proof of Proposition \ref{prop:0opt_is_MIGNE}}
\label{app:pf:prop:0opt_is_MIGNE}

Since $\bm u^\circ$ is a solution to the convexified problem in \eqref{eq:opt_main_cvx}, $\bm u^\circ$ satisfies all the constraints of Problem \eqref{eq:opt_main}, except possibly the binary constraints \eqref{eq:bin_const}. Thus, the tuple $(\{(p_i^{\mathrm{dg}}, d_i^{\mathrm{gu}})^\circ\}_{i \in \mc I^{\mathrm{dg}}}, \{\theta_i^\circ\}_{i \in \mc B}, \{(g_i^{\mathrm{s}})^\circ\}_{i \in \mc I^{\mathrm{gs}}},   \{\phi_{(i,j)}^\circ\}_{(i,j) \in \mc P^{\mathrm t}})$, which does not change after the second stage, satisfies \eqref{eq:p_pu_bound}, \eqref{eq:conv_dgu_ppu}, \eqref{eq:pow_bal}--\eqref{eq:gs_bound}, \eqref{eq:gbal},  \eqref{eq:gf_recip} and \eqref{eq:phi_lim}. Since $\tilde{\bm z}$ satisfies \eqref{eq:solve_binary} while $(\tilde y_{(i,j)}^{\psi_i}, \{\tilde y_{(i,j)}^m \}_{m=1}^r)$ satisfies \eqref{eq:y_aux_update}, the only constraint in the mixed-integer linear gas flow model that may not be satisfied is \eqref{eq:gf_eq6}. Therefore, if the recomputed pressure variable $\tilde{\bm \psi}$, which is a solution to Problem \eqref{eq:qp_pressure},  has zero optimal value, i.e., $J_\psi(\tilde{\bm \psi})=0$, then the tuple $(\phi_{(i,j)}^\circ, \tilde \psi_i, \tilde \psi_j, \tilde z_{(i,j)})$, for each  $(i,j) \in \mc P^{\mathrm{nt}}$, satisfies  \eqref{eq:gf_eq6} (and \eqref{eq:pres_lim}), thus implying that $\bm u^\star$ is a feasible point of Problem \eqref{eq:opt_main}. Since Problem \eqref{eq:opt_main_cvx} is a convex relaxation of Problem \eqref{eq:opt_main}, the value of the cost function in \eqref{eq:cost_f} evaluated at $\bm u^\circ$ is a lower bound of the optimal value of Problem \eqref{eq:opt_main}. Furthermore, since the decision variables that influence the cost functions, i.e., $((\{(p_i^{\mathrm{dg}})^\circ\}_{i \in \mc I^{\mathrm{dg}}},\{(g_i^{\mathrm{s}})^\circ\}_{i \in \mc I^{\mathrm{gs}}})$, are computed in the first stage and remain the same after the second stage, the value of the cost function in \eqref{eq:cost_f} after the second stage is the same as that after the first stage. Consequently, the lower bound is achieved by $\bm u^\star$. Thus, we conclude that $\bm u^\star$ is a solution to Problem \eqref{eq:opt_main}.

\bibliographystyle{IEEEtran}
\bibliography{ref}

\end{document}